\documentclass[10pt, a4paper]{article}
\usepackage{newtxtext,newtxmath}
\usepackage{amsmath,amsfonts}
\usepackage[ruled,vlined,linesnumbered]{algorithm2e}
\usepackage{array}
\usepackage[caption=false,font=normalsize,labelfont=sf,textfont=sf]{subfig}
\usepackage{textcomp}
\usepackage{stfloats}
\usepackage{url}
\usepackage{verbatim}
\usepackage{graphicx}
\usepackage{cite}
\usepackage{booktabs, tabularx}
\usepackage{eurosym}
\usepackage[nocomma]{optidef}

\newcommand{\variablestyle}[1]{\ensuremath{#1}}
\newcommand{\setstyle}[1]{\ensuremath{\mathcal{#1}}}
\newcommand{\numsetstyle}[1]{\ensuremath{\mathbb{#1}}}
\newcommand{\descriptorstyle}[1]{\ensuremath{\text{#1}}}
\newcommand{\vardescstyle}[2]{\variablestyle{#1}^{\descriptorstyle{#2}}}

\def\N{\setstyle{N}}
\def\A{\setstyle{A}}
\def\Ainv{\setstyle{A}^{\text{inv}}}
\def\Aninv{\setstyle{A}^{\text{not}}}
\def\Ain{\setstyle{A}^{\text{in}}}
\def\Aout{\setstyle{A}^{\text{out}}}
\def\G{\setstyle{G}}
\def\S{\setstyle{S}}
\def\SSeasonal{\setstyle{S}^{\text{seas}}}
\def\SShort{\setstyle{S}^{\text{ST}}}
\def\C{\setstyle{C}}
\def\D{\setstyle{D}}
\def\H{\setstyle{H}}

\def\L{\setstyle{L}}
\def\X{\setstyle{X}}
\def\Features{\setstyle{K}}
\newcommand{\Lfrom}[1]{\L^{\text{out}}_{#1}}
\newcommand{\Lto}[1]{\L^{\text{in}}_{#1}}
\def\R{\setstyle{R}}
\def\Reals{\numsetstyle{R}}

\def\vInvCost{\vardescstyle{c}{inv}}
\def\vOpeCost{\vardescstyle{c}{op}}
\def\vGenInv{\variablestyle{i}}
\def\vAcumUnit{\variablestyle{\pi}}
\def\vGenProd{\vardescstyle{p}{out}}
\def\vGenChar{\vardescstyle{p}{in}}

\def\vStorIntra{\vardescstyle{\sigma}{intra}}
\def\vStorIntraZero{\vardescstyle{\sigma}{intra,\,0}}
\def\vStorInter{\vardescstyle{\sigma}{inter}}
\def\vStorInterZero{\vardescstyle{\sigma}{inter,\,0}}
\def\vLineFlow{\variablestyle{f}}
\def\vSpillage{\vardescstyle{e}{sp}}
\def\vBorrow{\vardescstyle{e}{br}}

\newcommand{\pWeightRepr}[1]{\variablestyle{W}_{\!#1}}
\def\pWeight{\vardescstyle{W}{op}}
\def\pDemand{\variablestyle{D}}
\def\pDemandPeak{\vardescstyle{D}{max}}
\def\pGenAva{\variablestyle{A}}
\def\pInvCost{\variablestyle{I}}
\def\pVarCost{\variablestyle{V}}
\def\pVarCostSpillage{\vardescstyle{V}{sp}}
\def\pVarCostBorrow{\vardescstyle{V}{br}}
\def\pUnitCap{\variablestyle{P}}
\def\pIniUnit{\vardescstyle{U}{0}}
\def\pStorCap{\variablestyle{S}}
\def\pInflows{\variablestyle{E}}
\def\pInflowsPeak{\vardescstyle{E}{max}}
\def\pExpCap{\vardescstyle{L}{\!exp}}
\def\pImpCap{\vardescstyle{L}{\!imp}}
\def\pRamp{\variablestyle{R}}
\def\pEffIn{\vardescstyle{\eta}{in}}
\def\pEffOut{\vardescstyle{\eta}{out}}
\def\pEffOut{\vardescstyle{\eta}{out}}
\def\pInitialStorageLevel{\vardescstyle{S}{0}}
\def\numRPs{\vardescstyle{N}{RP}}

\newcommand{\setsize}[1]{\lvert#1\rvert}
\def\nD{\setsize{\setstyle{D}}}
\def\nH{\setsize{\setstyle{H}}}
\def\nR{\setsize{\setstyle{R}}}
\def\nT{\setsize{\setstyle{T}}}
\def\nY{\setsize{\setstyle{Y}}}

\def\placeholder{\,\cdot\,}
\def\clusteringMap{\variablestyle{\rho}}
\def\maxdist{\vardescstyle{\Delta}{\!max}}
\def\curdist{\vardescstyle{\Delta}{\!cur}}
\def\maxr{\vardescstyle{d}{max}}
\def\newr{\vardescstyle{r}{new}}

\def\ss{\vardescstyle{s}{se}}
\def\sn{\vardescstyle{s}{ST}}
\def\ainv{\vardescstyle{a}{inv}}
\def\aninv{\vardescstyle{a}{not}}

\def\feature{\variablestyle{k}}

\DeclareMathOperator{\proj}{proj}
\DeclareMathOperator{\pgd}{PGD}
\DeclareMathOperator{\dist}{dist}
\DeclareMathOperator{\argmax}{arg\,max}

\SetKw{Break}{break}

\begin{document}

\title{Hull Clustering with Blended Representative Periods for Energy System Optimization Models}

\author{%
Grigory Neustroev,
Diego A. Tejada-Arango,\\
Germ\'an Morales-Espa\~na,
Mathijs M. de Weerdt%
}

\maketitle

\begin{abstract}
The growing integration of renewable energy sources into power systems requires planning models to account for not only demand variability but also fluctuations in renewable availability during operational periods.
Capturing this temporal detail over long planning horizons can be computationally demanding or even intractable.
A common approach to address this challenge is to approximate the problem using a reduced set of selected time periods, known as representative periods (RPs).
However, using too few RPs can significantly degrade solution quality.
In this paper, we propose the method of hull clustering with blended RPs to enhance traditional clustering-based RP approaches in two key ways.
First, instead of selecting typical cluster centers (e.g., centroids or medoids) as RPs, our method is based on extreme points, which are more likely to be constraint-binding.
Second, it represents base periods as weighted combinations of RPs (e.g., convex or conic blends), approximating the full time horizon more accurately and with fewer RPs.
Through two case studies based on data from the European network operators, we demonstrate that hull clustering with blended RPs outperforms traditional RP techniques in both regret and computational efficiency.
\end{abstract}

\section{Introduction}

Energy System Optimization Models (ESOMs) are widely used to plan resilient and cost-effective energy systems.
Common examples in the power sector include generation expansion planning (GEP) \cite{Wogrin2022lego}, integrated resource planning \cite{IHLEMANN2022}, mid-term operational planning \cite{hellemo2024}, and transmission expansion planning \cite{RAMOS2022}.
As renewable energy sources like wind and solar gain traction, the variability and uncertainty in generation have increased---and will continue to rise with further integration \cite{IEA2021}.
Accurately capturing this variability requires high-resolution temporal modeling of demand and production.
However, incorporating such fine-grained data over long planning horizons greatly increases the computational burden, making ESOMs difficult---or even intractable---to solve for large-scale systems \cite{CHO2022}.

\subsection{Motivation}

To address this challenge, researchers and practitioners often apply reduction techniques to temporal data.
Instead of using all time series, they work with smaller subsets of time-varying data such as demand and renewable availability profiles.
When periods are sufficiently similar, this can significantly reduce runtime without introducing substantial errors.
Common approaches include time blocks and time slices \cite{Wogrin2014}, as well as representative periods (RPs)  \cite{Nahmmacher2016}.
Recently, the latter has gained popularity, as it keeps the high temporal resolution while preserving the chronological structure within RPs, allowing for the modeling of seasonal storage, inter-period ramping, and other temporal constraints \cite{TejadaArango2018,Kotzur2018b}.

The concept of RPs refers to specific time intervals, such as days or weeks, selected to capture the variability and characteristics of the ESOM over an extended period, such as a year or even decades.
Selecting RPs involves a trade-off between model accuracy and computational efficiency.
Therefore, identifying critical periods that accurately capture system behavior is an important research topic, as these periods influence how well the system is represented while still enabling accurate analysis \cite{Wogrin2019}.

\subsection{Literature Review}

Clustering algorithms are commonly utilized to obtain RPs based on time series of demand and renewable availability profiles \cite{Poncelet2017}.
Traditional clustering methods, such as \(k\)-medoids \cite{TejadaArango2018} and hierarchical clustering \cite{Pineda2018}, are among the most common and straightforward approaches in the literature.
However, these methods may not adequately capture extreme events or rare occurrences crucial for energy systems investment and operational planning decisions \cite{TEICHGRAEBER2020}.

To address this limitation, researchers have proposed assigning greater weight to extreme values in the dataset \cite{GarciaCerezo2022}.
This ensures that boundary conditions and rare events are captured in the reduced set of RPs.
These weights are typically defined a priori, either through heuristic rules or expert knowledge of the system.
Other approaches are data-driven.
For example, \cite{Sun2019} introduces a framework that selects RPs based on investment decisions rather than input data. 
However, this method requires solving multiple reduced model instances, which can be prohibitive---particularly for large-scale ESOMs with seasonal storage.
Alternatively, hybrid clustering techniques, such as combining \(k\)-means with hierarchical clustering, have also been explored to improve the selection of RPs \cite{Liu2018}.

The concept of RP mixing has emerged as a means to recover information from the non-represented periods as a linear combination of the RPs.
This concept was studied in \cite{Gonzato2021} to better capture the diversity of original periods to reflect a broader range of conditions in the RP version of the ESOM.
Nevertheless, the proposal in \cite{Gonzato2021} overlooked the impact of the extreme conditions in the representative selection and proposed a mixed-integer programming optimization to find the RPs, which can be computationally intensive for large-scale problems.
Additionally, that approach does not take into account the fact that the RPs will be mixed with each other in the expansion planning problem.

In conclusion, challenges persist in developing temporal aggregation methods that can effectively reduce computational requirements without significant loss of accuracy or that depend either on specific predefined rules or previous solutions of the reduced versions of the ESOM.
The need for approaches that can capture both typical and extreme system conditions computationally efficiently is increasingly recognized.

\subsection{Contributions}

In response to the challenges shown in the previous section, this paper proposes the \emph{hull clustering with blended representative periods} that aims to capture extreme conditions in the input data with fewer representatives than the counterparts, efficiently achieving high levels of accuracy, while maintaining a high fidelity when representing the non-extreme periods.
The fundamental elements of the hull clustering with blended RPs are as follows.
\begin{enumerate}
  \item \emph{Selection of Extreme Points:} Instead of using cluster centroids or medoids, blended RPs select representatives from the data hull, focusing on extreme data points. This ensures that critical scenarios, which are often constraint-binding, are adequately represented, leading to optimization outcomes with lower regret.
  \item \emph{Conic Sums for Blending:} By using conic sums of representative periods, blended RPs create weighted combinations that more accurately reflect the data in the base periods. This approach reduces the number of RPs needed while preserving the integrity of the original data distribution, thus enhancing computational efficiency without compromising accuracy.
\end{enumerate}

The positive effects of the contributions of this paper include:
\begin{enumerate}
  \item \emph{Improved Solution Quality:} The use of extreme points ensures that the selected RPs better capture the variability of the data, leading to more accurate solutions.
  \item \emph{Reduction in Computational Burden:} By reducing the number of RPs needed, the method lowers computational requirements, making longer planning horizons and more detailed models feasible.
\end{enumerate}
However, the adoption of blended RPs introduces complexities, as the need for extreme points challenges conventional clustering techniques and demands the development of specialized algorithms to implement the hull clustering effectively.

In addition, we conduct two case studies, with data from the European network.
The first one is a stylized version of a generation expansion planning problem, considering ramping but not storage.
The second is a Power-to-X dispatch case study based on the global ambition scenario of the TYNDP for the year 2050 \cite{tyndp}, which includes seasonal storage for hydro reservoir technologies.
We demonstrate that our method achieves the same solution quality as conventional methods, using substantially fewer RPs and therefore runtime.
\section{Model Formulation}

Although the hull clustering with the blended RP method can be applied to a wide range of ESOMs, this paper focuses on two versions:
a generation expansion planning (GEP) and
a sector-coupling power and hydrogen dispatch problem (Power-to-X, P2X). 
Both versions include typical complexities in energy systems, such as transmission capacities, energy conversion, short-term and seasonal storage technologies, and ramping constraints, to illustrate the proposed method's core ideas.
In this section, we present the general formulation of an ESOM that can be adapted for both GEP and P2X optimization models.
The GEP model aims to determine the optimal investment in new generation capacity while minimizing the sum of investment and operational costs, subject to changing demand, ramping, and resource availability.
The P2X dispatch model aims to determine the optimal dispatch of existing capacities while minimizing operational costs, subject to available resource capacities, balance, conversion, and storage constraints.
Given the complexity and temporal resolution of these problems, reducing the time series data through RPs is a common approach.
Here, we outline the ESOM formulation and how RP methods are used to make large-scale problems computationally tractable.

\subsection{Energy System Optimization Model}\label{sec:ESOM}

In energy system optimization, we consider a set of assets \(a \in \A\).
Each asset is located at a node \(n \in \N\) in a network, which is connected by transmission lines \(l \in \L\) transporting energy carriers \(x \in \X\) (e.g., electricity or gas).
We use \(\Ain_{n,x}\) and \(\Aout_{n,x}\) for assets that respectively consume and output carrier \(x\) at node \(n\).
These assets are categorized into three functional types:
\begin{enumerate}
  \item generation assets \(g \in \G\), which produce energy (e.g., wind turbines, solar panels);
  \item storage assets \(s \in \S\), which store excess energy for later use (e.g., batteries, pumped hydro), divided into short-term \(\sn \in \SShort\) and seasonal \(\ss \in \SSeasonal\) storage; and
  \item conversion assets \(c \in \C\), which transform energy between carriers (e.g., fuel cells).
\end{enumerate}
The asset types are disjoint:
\[\A = \G \sqcup \S \sqcup \C,\quad\text{where}\quad \S = \SShort \sqcup \SSeasonal.\]
Some assets \(\ainv \in \Ainv \subseteq A\) are designated as investable, meaning their capacities can be expanded as part of the planning process, subject to investment costs and budget constraints.
Non-investable assets \(\aninv \in \Aninv = A \setminus \Ainv\), by contrast, represent existing infrastructure. In addition, we define \(\Lfrom{n,x}\) and \(\Lto{n,x}\) to indicate the subset of outgoing and incoming lines to the node and carrier, respectively.

To capture the system's temporal dynamics, optimization is performed over a horizon \(\nT = \nD \cdot \nH\), which is discretized into a set of \(\nD\) time periods \(d \in \D\) (e.g., days) of equal length \(\nH\), each consisting of time steps \(h \in \H\) (e.g., hours).

The temporal dimension of the model can be reduced by replacing base periods \(d \in \D\) with a much smaller set of RPs \(r \in \R\), where \(\nR \ll \nD\).
The RPs are selected so that they reflect typical combinations of the relevant conditions, reducing the computational burden while preserving the essential temporal patterns needed for accurate decision-making.

RPs are typically identified using clustering methods.
While \(k\)-means is commonly used, alternative methods such as \(k\)-medoids and hierarchical clustering may offer better performance depending on the structure of the data.
For instance, \(k\)-medoids selects actual observations, making it more robust to outliers.
After clustering, each base period \(d\) is assigned an RP \(d \mapsto \clusteringMap(d)\) by a clustering function \(\clusteringMap: \D \to \R\).
As a result, they are replaced with their RPs \(r = \clusteringMap(d)\).
We introduce weights \(\pWeightRepr{d,r}\) as indicators of this assignment:
\begin{equation}
  \text{if } \clusteringMap(d) = r \text{ then } \pWeightRepr{d,r} = 1, \text{ otherwise } \pWeightRepr{d,r} = 0.
  \label{eq:dirac-weights}
\end{equation}
Additionally, RPs \(r\) are assigned total weights \(\pWeightRepr{r}\) representing the number of periods that are replaced by this representative:
After clustering, each base period \(d\) is assigned a representative $r$ 
%\(d \mapsto \clusteringMap(d)\) 
by a clustering function \(\clusteringMap: \D \to \R\). 
In the model, each base period \(d\) is then replaced by its representative \(r = \clusteringMap(d)\).
RPs \(r\) are assigned weights \(\pWeightRepr{r} = \sum_{d \in \D} \pWeightRepr{d,r}\) representing the number of periods that are replaced by this representative.

Using this approach, we propose the following ESOM that integrates investment and operational decisions.
All input parameters are denoted by capital letters and decision variables by lower case letters.
Except for the flows \(\vLineFlow_{l, r, h}\), the decision variables are constrained to be non-negative, but these constraints are omitted for brevity.

\begin{mini!}[3]<b>
{}%
{\vInvCost +\vOpeCost\label{obj:esom}}%
{\label{esom}}%
{}
% constraints
%total investment cost
\addConstraint%
  {\vInvCost}%
  {= \smashoperator{\sum_{\ainv \in \Ainv}} \pInvCost_{\ainv} \pUnitCap_{\ainv} \vGenInv_{\ainv}\label{cnst:investment-cost}}%
  {}
% total operation cost
\addConstraint%
  {\vOpeCost}%
  {= \pWeight \sum_{r \in \R} \pWeightRepr{r} \sum_{h \in \H} \Bigl( %
    \sum_{g \in \G} \pVarCost_{g} \vGenProd_{g, r, h}%
    + \makebox[0pt][l]{$\tfrac{1}{\tau}\smashoperator{\sum_{\ss \in \SSeasonal}} \bigl(\pVarCostSpillage_{\ss} \vSpillage_{\ss, r, h}
    + \pVarCostBorrow_{\ss} \vBorrow_{\ss, r, h}\bigr)\Bigr)$}%
  }%
  {\label{cnst:operation-cost}}
% node balance
\addConstraint%
  {%
      \smashoperator{\sum_{a \in \Ain_{n,x}}} \vGenProd_{a, r, h}%
    - \smashoperator{\sum_{a \in \Aout_{n,x}}} \vGenChar_{a, r, h}%
    - \smashoperator{\sum_{l \in \Lfrom{n,x}}}\vLineFlow_{l, r, h} + \smashoperator{\sum_{l \in \Lto{n,x}}} \vLineFlow_{l, r, h}%
    \label{cnst:node-balance}
  }
  {
    = \pDemand_{n, x, r, h}\pDemandPeak_{n, x}\quad
  }
  {\forall n, x, r, h,}
% intra-storage balance
\addConstraint%
  { \vStorIntra_{\sn, r, 1} - \vStorIntraZero_{\sn, r}}%
  {= \bigl(\pEffIn_{\sn} \vGenChar_{\sn, r, 1}%
    - \tfrac{\vGenProd_{\sn, r, 1}}{\pEffOut_{\sn}}\bigr)\tau%
    \label{cnst:intra-short-term-storage-balance-0}
  }
  {\forall \sn, r,}
\addConstraint%
  { \vStorIntra_{\sn, r, h} - \vStorIntra_{\sn, r, h-1}}%
  {= \bigl(\pEffIn_{\sn} \vGenChar_{\sn, r, h}%
    - \tfrac{\vGenProd_{\sn, r, h}}{\pEffOut_{\sn}}\bigr) \tau%
    \label{cnst:intra-short-term-storage-balance}}
  {\forall \sn, r, h \neq 1,}
\addConstraint%
  { \vStorIntra_{\ss, r, 1} - \vStorIntraZero_{\ss, r}}%
  {= \makebox[0pt][l]{$\bigl(\pEffIn_{s} \vGenChar_{\ss, r, 1}%
    - \frac{\vGenProd_{\ss, r, 1}}{\pEffOut_{\ss}}\bigr)\tau%
    -\vSpillage_{\ss, r, 1} + \vBorrow_{\ss, r, 1}+ \pInflows_{\ss, r, 1}\pInflowsPeak_{\ss}$}%
    \nonumber
  }
  {}
\addConstraint%
  {}
  {\label{cnst:intra-seasonal-storage-balance-0}}
  {\forall \ss, r,}
\addConstraint%
  { \vStorIntra_{\ss, r, h} - \vStorIntra_{\ss, r, h-1}}%
  {= \makebox[0pt][l]{$\bigl(\pEffIn_{\ss} \vGenChar_{\ss, r, h}%
    - \frac{\vGenProd_{\ss, r, h}}{\pEffOut_{\ss}}\bigr)\tau -\vSpillage_{\ss, r, h} +\vBorrow_{\ss, r, h} + \pInflows_{\ss, r, h}\pInflowsPeak_{\ss}$}%
    \nonumber
  }
  {}
\addConstraint%
  {}
  {\label{cnst:intra-seasonal-storage-balance}}
  {\forall \ss, r, h \neq 1,}
%% inter-storage balance
\addConstraint%
  { \vStorInter_{\ss, 1} - \vStorInterZero_{\ss}}%
  {= \smashoperator{\sum_{r \in \R}} \pWeightRepr{1,r} (\vStorIntra_{\ss, r, H} - \vStorIntraZero_{\ss, r})%
    \label{cnst:inter-seasonal-storage-balance-0}%
  }
  {\forall \ss}
\addConstraint%
  { \vStorInter_{\ss, d} - \vStorInter_{\ss, d-1}}%
  {= \smashoperator{\sum_{r \in \R}} \pWeightRepr{d,r} (\vStorIntra_{\ss, r, H} - \vStorIntraZero_{\ss, r})%
    \label{cnst:inter-seasonal-storage-balance}%
  }
  {\forall \ss, d \neq 1,}
%% storage cycling
\addConstraint%
  { \vStorInterZero_{\ss}}%
  {%
   = \pInitialStorageLevel_{\ss}%
   = \vStorInter_{\ss, D}%
   = \smashoperator{\sum_{r \in \R}} \pWeightRepr{D,r}\vStorIntra_{\ss, r, H}
   \label{cnst:inter-storage-cycling}%
  }
  {\forall \ss,}
\addConstraint%
  { \vStorIntra_{\sn, r, H}}%
  {= \vStorIntraZero_{\sn, r}%
  }
  {\forall \sn, r%
    \label{cnst:intra-storage-cycling}} 
% conversion balance
\addConstraint%
  { \pEffIn_{c} \vGenChar_{c, r, h}}%
  {= \tfrac{\vGenProd_{c, r, h}}{\pEffOut_{c}}}%
  {\forall c, r, h,%
    \label{cnst:conversion-balance}
  }
% accumulated units investable
\addConstraint%
  {\vAcumUnit_{\ainv}}%
  {=  \pUnitCap_{\ainv} (\pIniUnit_{\ainv} + \vGenInv_{\ainv})}%
  {%
    \forall \ainv,%
    \label{cnst:acummulated-units-investable}%
  }
% accumulated units not investable
\addConstraint%
  {\vAcumUnit_{\aninv}}%
  {= \pUnitCap_{\aninv} \pIniUnit_{\aninv}}%
  {%
    \forall \aninv,%
    \label{cnst:acummulated-units-non-investable}%
  }
% maximum production
\addConstraint%
  {%0 \leq
    \vGenProd_{a, r, h}%
  }%
  {\leq \pGenAva_{a, r, h} \vAcumUnit_{a}}%
  {%
    \forall a, r, h,%
    \label{cnst:max-production}%
  }
% maximum consumption
\addConstraint%
  {%0 \leq
    \vGenChar_{s, r, h}%
  }%
  {\leq \vAcumUnit_{s}}%
  {%
    \forall s, r, h,%
    \label{cnst:max-consumption}%
  }
% ramping intra
\addConstraint%
  {%
    -\pRamp_{g} \vAcumUnit_{g}\tau \leq \vGenProd_{g, r, h} - \vGenProd_{g, r, h-1}%
  }%
  {\leq \pRamp_{g} \vAcumUnit_{g} \tau \quad}%
  {\forall g, r, h \neq 1,%
    \label{cnst:ramp-intra}%
  }
% ramping inter
\addConstraint%
  {-\pRamp_{g} \vAcumUnit_{g}\tau \leq \smashoperator{\sum_{r \in \R}} (\pWeightRepr{d,r}\vGenProd_{g, r, 1} - \pWeightRepr{d-1,r}\vGenProd_{g, r, H})%
  }
  {%
    \leq \pRamp_{g} \vAcumUnit_{g}\tau\label{cnst:ramp-inter}%
  }
  {\forall g, d \neq 1,}
% domains
\addConstraint%
  {%0 \leq
    \vStorIntra_{s, r, h}%
  }%
  {\leq \pStorCap_{s}}%
  {\forall s, r, h,%
    \label{cnst:intra-max-energy-storage}%
  }
\addConstraint%
  {%
    \pStorCap_{\ss, d}^{\text{min}} \pStorCap_{\ss} 	\leq
    \vStorInter_{\ss, d}%
  }%
  {\leq \pStorCap_{\ss, d}^{\text{max}} \pStorCap_{\ss}}%
  {\forall \ss, d,%
    \label{cnst:inter-max-energy-storage}%
  }
\addConstraint%
  {%
    -\pImpCap_{l} \leq \vLineFlow_{l, r, h}%
  }%
  {\leq \pExpCap_{l}}%
  {%
    \forall l, r, h.%
    \label{cnst:flow-domain}%
  }
\end{mini!}

The model \eqref{esom} minimizes the total cost \eqref{obj:esom} comprising two main components.
First, the \emph{total investment cost} \(\vInvCost\) in \eqref{cnst:investment-cost} is computed as the sum of investment costs across all investable assets \(a \in \Ainv\).
For each asset, the cost is determined by its capacity per unit \(\pUnitCap_{a}\), the annualized investment cost per MW \(\pInvCost_{a}\), and the number of units selected for investment \(\vGenInv_{a}\).
Second, the \emph{total operational cost} \(\vOpeCost\) in \eqref{cnst:operation-cost} accounts for generation and reliability-related costs across all RPs \(r \in \R\) and time steps \(h \in \H\).
These include costs associated with generation output \(\vGenProd_{g,r,h}\), as well as penalties for spillage \(\vSpillage_{\ss,r,h}\) and borrowing \(\vBorrow_{\ss,r,h}\) of water.
Spillage and borrowing are slack variables for water reservoir level constraints; for non-reservoir technologies, these are set to zero.
To ensure consistency with annualized investment cost, the operational cost is also annualized using a weight \(\pWeight = \nY / \nT\), where \(\nY\) and \(\nT\) are the numbers of time steps in a year and in the model, and \(\tau\) is the duration of the time step.

The next group of constraints enforces intra-period energy system balances.
The \emph{node balance constraint} \eqref{cnst:node-balance} ensures that the local demand \(\pDemand_{n,x,r,h}\pDemandPeak_{n,x}\) is satisfied by the sum of locally generated, consumed, and net imported power (imports minus exports).
Here \(\pDemandPeak_{n,x}\) represents the peak demands, and \(\pDemand_{n,x,r,h} \in [0, 1]\) are the time-varying demand profiles.
The \emph{intra-period storage balance constraints} \eqref{cnst:intra-short-term-storage-balance-0}--\eqref{cnst:intra-seasonal-storage-balance} track the state of charge \(\vStorIntra_{s,r,h}\) of storage assets across time, accounting for the previous state of charge, initial storage level \( \vStorIntraZero_{s,r} \), charging and discharging flows adjusted by efficiencies \(\pEffIn_{s}\) and \(\pEffOut_{s}\), spillage \(\vSpillage_{s,r,h}\) and borrowing \(\vBorrow_{s,r,h}\), and external inflows \(\pInflows_{s,r,h}\pInflowsPeak_{s,r,h}\), where \(\pInflows_{s,r,h} \in [0,1]\), commonly representing hydro inflows for reservoirs. 
Some technologies may not include either internal \(\pEffIn_{s} \vGenChar_{s, r, h} \tau\) or external inflows, in which case either \(\vGenChar_{s, r, h}\) or \(\pInflowsPeak_{s,r,h}\) are set to zero.

The \emph{inter-period storage balance constraints} \eqref{cnst:inter-seasonal-storage-balance-0}--\eqref{cnst:inter-seasonal-storage-balance} recover the chronological information for the seasonal storage when using representative periods by linking the inter- and intra-period storage levels through the use of the weights \(\pWeightRepr{d,r} \) (note that the weights are not necessarily integer).

The \emph{cyclic constraints} \eqref{cnst:inter-storage-cycling} link the last inter-period storage level \(\vStorInter_{\ss,D}\) to the initial one \(\vStorInterZero_{\ss}\) and set both to a pre-defined value \(\pInitialStorageLevel_{\ss}\). This is done to avoid creating a solution where the storage is fully charged initially and discharges over the time horizon, creating free energy in the system. Intra-period cycling constraints \eqref{cnst:intra-storage-cycling} serve the same role but for short-term storage.
Additionally, in the model formulation so far, the intra-period storage variables \(\vStorIntra_{s,r,h}\) are defined up to a constant: increasing all or them by the same amount does not change the inter-period variables \(\vStorInter_{s,d}\). To address this issue, we ``tether'' them to each other via the last part of \eqref{cnst:inter-storage-cycling}.

The \emph{conversion balance constraint} \eqref{cnst:conversion-balance} ensures that energy inputs and outputs at conversion assets remain consistent with their conversion efficiencies \(\pEffIn_{c}\) and \(\pEffOut_{c}\).

The next set of constraints captures the \emph{physical limitations of the assets}.
Constraints \eqref{cnst:acummulated-units-investable} and \eqref{cnst:acummulated-units-non-investable} define the total number of available units \(\vAcumUnit_{a}\) as the sum of existing units \(\pIniUnit_{a}\) and newly invested units \(\vGenInv_{a}\).
Operational capacity constraints then limit each asset's activity based on this availability.
Specifically, constraints \eqref{cnst:max-production} and \eqref{cnst:max-consumption} ensure that the produced power \(\vGenProd_{a,r,h}\) and consumed power \(\vGenChar_{a,r,h}\) do not exceed the time-dependent availability \(\pGenAva_{a,r,h}\).
This parameter equals 1 for all assets except renewable generators, for which it reflects hourly availability with \(\pGenAva_{g,r,h}\in [0,1]\).
Ramping constraints \eqref{cnst:ramp-intra} and \eqref{cnst:ramp-inter} limit the change in generation between consecutive time steps by the ramping parameter \(\pRamp_{g}\), modeling real-world operational inflexibility.
Constraints \eqref{cnst:intra-max-energy-storage} and \eqref{cnst:inter-max-energy-storage} ensure that storage levels remain withing required upper and lower limits, while transmission constraints \eqref{cnst:flow-domain} enforce that line flows \(\vLineFlow_{l,r,h}\) stay within import \(-\pImpCap_{l}\) and export \(\pExpCap_{l}\) limits.

In this paper, we focus on this specific problem formulation for clarity of presentation.
When the investment decisions \(\vGenInv_{a}\) are treated as decision variables and restricted to generation assets---i.e., \(\Ainv \subseteq \G\)---the model \eqref{esom} reduces to the GEP problem.
Alternatively, if all investment decisions \(\vGenInv_{a}\) are fixed as input parameters, the model can be used for \emph{operations planning} such as P2X, focusing solely on the cost-efficient dispatch of existing generation assets.
The formulation is flexible and can be used for \emph{storage investment}, and with trivial adjustment extended to \emph{transmission expansion planning}.
\section{Methodology}

We have identified that any representative period method consists of two key components: a clustering method for finding RPs, and the mapping of base-period data to RPs.
In this paper, we introduce a generalized framework that allows for flexible choices in both components.
This leads to a new variant, which we term the \emph{method of blended RPs}.
The resulting procedure is summarized in Algorithm~\ref{alg:blended}.

\begin{algorithm}[tbh]
\caption{The Method of Blended RPs\label{alg:blended}}
\KwData{%
  time-varying clustering data,
  number of representative periods,
  clustering method,
  weight type.
}
\KwResult{approximate solution of ESOM \eqref{esom}.}
collect the clustering data into a clustering matrix \(\mathbf{C}\)\;
find the RP data \(\mathbf{R}\) via the chosen clustering method\;
find weights \(\mathbf{W}\) of the chosen type\;
solve the reduced ESOM using the RP data \(\mathbf{R}\)\;
\Return{the resulting solution.}
\end{algorithm}

\subsection{Overview of Blending in Representative Periods}

The first key innovation of the method of blended representative periods lies in generalizing how base periods are represented in terms of a reduced set of RPs.
Traditional clustering-based approaches assign each base period to exactly one RP.
In such a framework, the clustering function \(\clusteringMap: \D \to \R\) defines a partition of the original periods: each base period \(d \in \D\) is mapped to a single RP \(r \in \R\).
This assignment, however, can be limiting, as in many practical settings---such as renewable generation---the base periods exhibit patterns that are mixtures of different underlying conditions.
For instance, a mildly sunny, breezy afternoon may not be well captured by any single RP but may lie somewhere between a typical sunny day and a typical windy day.
In such cases, forcing a one-to-one mapping can introduce large approximation errors.

To address this, we introduce the idea of blended representative periods, where each base period is approximated by a weighted combination of multiple RPs.
This blend is captured by a weight matrix\footnote{In matrix definitions, we use subscript and superscript to describe row and column indices respectively.}
\(
  \mathbf{W} = [W_{d,r}]_{d \in \D}^{r \in \R},
\)
where each row specifies how the base period \(d\) is expressed in terms of RPs \(r\).
Instead of mapping each base period \(d\) to a single RP \(r\), we can now represent the problem data \(z\) as
\(
  y_d \approx \sum_{r \in \R} \pWeightRepr{d,r} \cdot y_r,
\)
where \(y_d\) stands for a period-dependent parameter or variable, such as demand \(\pDemand_{n, x, d, h},\) or produced power \(\vGenProd_{a, d, h}\).
This generalization allows us to interpolate between RPs and better preserve the diversity and nuance of the original time series.

\subsection{Blending Types: Dirac, Convex, and (Sub-Unit) Conical}

The choice of weight matrix \(\mathbf{W}\) is central to how base periods are approximated in the blended RP method.
We have identified and analyzed four types of blending schemes, based on the geometric structure imposed on the weights: Dirac, convex, sub-unit conical, and conical.
These represent increasing levels of flexibility in how base periods can relate to the representative periods.

\subsubsection{Dirac Weights (`Hard' Clustering)}
This is the traditional approach discussed in Section~\ref{sec:ESOM}.
Each base period \(d \in \D\) is assigned to exactly one RP \(r \in \R\), and the weight matrix \(\mathbf{W}\) is binary with a single nonzero entry per row as per \eqref{eq:dirac-weights}.
This structure corresponds to a Dirac measure in discrete space (i.e., a measure concentrated on a single point).
While computationally efficient, Dirac blending cannot capture overlap between RPs and may introduce abrupt transitions.

\subsubsection{Convex Weights}
A more flexible alternative allows each base period to be written as a \emph{convex combination} of RPs:
\[
  \sum_{r \in \R} \pWeightRepr{d,r} = 1, \quad \pWeightRepr{d,r} \geq 0, \quad\forall r \in \R, d \in \D.
\]
Convex blending allows smooth interpolation between RPs; as such, it has already been used in the context of ESOMs \cite{Gonzato2021}.

\subsubsection{Sub-Unit Conic Weights}\label{sec:sub-unit-conical}
This more expressive variant allows the weights to be non-negative and sum at most into 1:
\[
  \sum_{r \in \R} \pWeightRepr{d,r} \leq 1, \quad \pWeightRepr{d,r} \geq 0, \quad\forall r \in \R, d \in \D.
\]
This ensures that inequality constraints remain valid without modification.
To see why, consider an arbitrary linear inequality involving variables \(y_{i,r}\) associated with RPs,
\(
  \sum_{i} a_i y_{i,r} \leq b, \quad \forall r \in \R,
\)
where \(a_i\) are arbitrary constraint coefficients and \(b\) is non-negative.
This implies that \(\sum_{i} a_i y_{i,d} \leq b\) in the base problem holds as well:
\begin{align}
  \sum_{i} a_i y_{i,d}
  &= \sum_{i} a_i  \sum_{r \in \R} \pWeightRepr{d,r} y_{i,r}\nonumber\\
  &= \sum_{r \in \R} \pWeightRepr{d,r} \underbrace{\sum_{i} a_i y_{i,r}}_{\leq b} \leq b \underbrace{\sum_{r \in \R} \pWeightRepr{d,r}}_{\leq 1} \leq b.\label{eq:upper-bound-sub-unit}
\end{align}
By the same argument, inequalities of the form \(-b \leq \sum_{i} a_i y_{i,d}\) hold for non-negative constants \(b\).
Thus, sub-unit conic blending naturally preserves the feasibility of inequality constraints in the blended problem, \emph{without any additional reformulation}.

\subsubsection{Conic Weights}\label{sec:conical}
The most general linear blending we consider imposes only non-negativity:
\[
  \pWeightRepr{d,r} \geq 0, \quad\forall r \in \R, d \in \D.
\]
This allows each base period to be approximated by any non-negative linear combination of RPs, including combinations with total weight greater than one.
Such conic blending is the most expressive but may be less interpretable or lead to infeasibility of the solution to the original problem with respect to upper bounds, as the reasoning used in \eqref{eq:upper-bound-sub-unit} no longer holds.

\medskip

As far as we are aware, the (sub-unit) conic weights have not been described in the literature in the context of RPs.

Geometrically, each type of blending weight corresponds to projecting base period data onto a different type of hull formed by the RPs.
This results in projection errors, which are the differences between the original base period data and its reconstruction from a weighted combination of representative periods, and it decreases as the space of admissible weights becomes more expressive.
Specifically, Dirac weights select a single RP corresponding to the vertices of the set.
Convex weights project each base period into the convex hull of the RPs; data points inside the hull can be approximated without projection error.
Sub-unit conic weights extend the hull by adding the null point into it, reducing projection errors even further.
Finally, general conic weights project base periods into the full conic hull of the RPs.

\begin{figure*}[tbh]
\centering
\includegraphics[width=\textwidth]{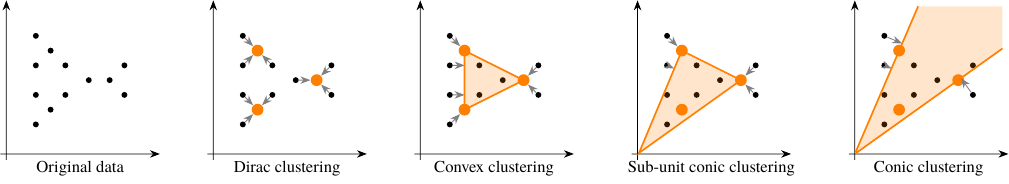}
\caption{
    Projection errors when approximating base period data using different weight types.
    Orange areas show the spaces of all points which can be represented without introducing a projection error.
    Errors decrease from left to right as we move from a discrete Dirac to more general weight types.
  }
  \label{fig:projection-error}
\end{figure*}

Figure~\ref{fig:projection-error} illustrates this visually.
It shows how the projection error shrinks as we move from approximating a base period with a single RP (Dirac) to the full conic hull.
The more expressive the hull, the more accurately we can approximate the data that is being clustered.
This hierarchy reflects a trade-off between expressiveness and constraint preservation.

We use the term \emph{blended} weights as a general descriptor encompassing all four weight types introduced here: Dirac, convex, sub-unit conic, and conic weights. These sets of weights satisfy the strict inclusions
\[
  \mathbb{W}^{\text{Dirac}} \subset \mathbb{W}^{\text{conv}} \subset \mathbb{W}^{\leq 1} \subset \mathbb{W},
\]
where \(\mathbb{W}^{\placeholder}\) is the space of all weight matrices \(\mathbf{W}\) of a type \((\placeholder)\) and \(\mathbb{W}\) is the space of conic weights.
Thus, each subsequent class further generalizes the conventional Dirac weights.

Some of these spaces have well-known geometric interpretations and standard names. For instance, the space of all conic weights corresponds to the positive orthant \(\Reals_+^n\), \(\mathbb{W} = \Reals_+^n\), the space of convex weights is the standard \(n\)-simplex \(\bigtriangleup_n\), \(\mathbb{W}^{\text{conv}} = \bigtriangleup_n = \bigl\{\mathbf{v}\in\Reals^n_+ \bigm| \lVert\mathbf{v}\rVert_1 = 1\bigr\}\), and the space of Dirac weights is the standard \(n\)-dimensional basis \(\mathbb{W}^{\text{Dirac}} = \{\mathbf{e}_i\}_{i=1, \dots, n}\).

Besides providing a unifying framework that generalizes existing RP methods, we propose to blend these periods using weights that are not necessarily convex.
This allows our method to go beyond simple interpolation between RPs and instead express the base period data with greater flexibility and nuance.

\subsection{Constructing the Clustering Matrix}

The first step of Algorithm~\ref{alg:blended} is to represent each base period \(d \in \D\) by a feature vector that captures its relevant time-varying characteristics.
Collecting these vectors forms the clustering matrix \(\mathbf{C}\).
This matrix is required because clustering algorithms such as \(k\)-medoids operate in a feature space:
to meaningfully group periods, each period must be described by a vector that reflects its time-varying data, rather than by its index alone.
For ESOMs, these typically include demands \(\pDemand_{n,x,d,h}\), availability profiles \(\pGenAva_{g,d,h}\) of renewable generation technologies, and time-varying energy storage inflows \(\pInflows_{s,d,h}\).

Formally, we define
\(
  \mathbf{C} = [c_{\feature,d}]_{\feature \in \Features}^{d \in \D}
\)
where each column \(\mathbf{c}_{d}\) encodes the time series of all features for period \(d\).
To obtain these feature vectors, we first compute the data matrices for each feature type and then stack them vertically:
\begin{equation*}
\left.
\begin{aligned}[c]
  \mathbf{D} &= [\pDemand_{n,x,d,h}]_{(n,x,h)\in \N \times \X \times \H}^{d \in \D},\\
  \mathbf{A} &= [\pGenAva_{g,d,h}]_{(g,h)\in \G \times \H}^{d \in \D},\\
  \mathbf{E} &= [\pInflows_{s,d,h}]_{(s,h)\in \S \times \H}^{d \in \D},\\
\end{aligned}
\quad\right\}\,\longrightarrow\quad
\begin{aligned}[c]
\mathbf{C} = \left[\begin{array}{c}
  \mathbf{D} \\
  \hline 
  \mathbf{A} \\
  \hline
  \mathbf{E}
\end{array}\right].
\end{aligned}
\end{equation*}
To ensure that all features contribute comparably to the clustering process, we use normalized data;
that is, the problem is formulated such that all values used in clustering lie between zero and one.
Once \(\mathbf{C}\) is constructed, it is passed to the clustering algorithm of choice to find the RPs.

\subsection{Hull Clustering}
Conventional clustering algorithms such as \(k\)-means and \(k\)-medoids are typically well-suited for Dirac-weight approximations, where each base period is assigned to a single RP.
These methods minimize average point-to-cluster-center distances but do not consider whether the selected cluster centers (RPs) span a suitable hull for blending.

When using \emph{blended weights}, it becomes essential that each base period can be reconstructed as a weighted combination of RPs, where the weights satisfy specific constraints.
If the selected RPs do not span the correct hull, base periods will be projected outside of the feasible blending space, resulting in significant projection error.

\subsubsection{Hull Types}

Given a set of RPs, the choice of weight type determines the kind of hull within which the original periods must lie in order to be exactly reconstructible. Specifically:

\begin{itemize}
  \item \emph{Convex hull}: For convex weights, the original periods must lie within the convex hull of the representative periods. This is the smallest convex set containing them.
  
  \item \emph{Convex hull with null}: For sub-unit conic weights, reconstruction is possible if the original periods lie within the convex hull of the representative periods and the origin (null). Intuitively, this corresponds to finding convex weights that may assign some mass to a null (zero) vector. After discarding the null, the remaining positive weights sum to less than one, hence sub-unit.
  
  \item \emph{Conic hull}: For general conic weights, the reconstruction is valid as long as the original periods lie within the conic hull of the RPs.
\end{itemize}

\begin{figure*}[tbh]
\centering
\includegraphics[width=\textwidth]{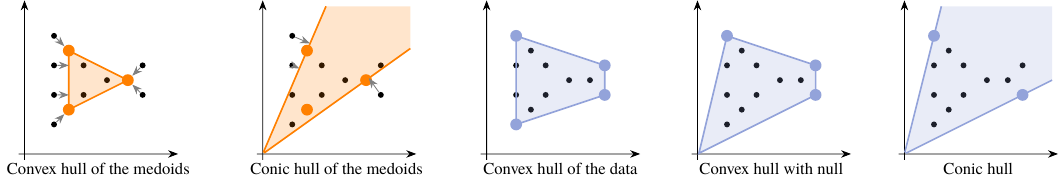}
\caption{
    Geometric interpretation of different hull types. 
    A set of base period data (dots) is shown, and the shaded region indicates the span of each hull.
    Even when blended weights are used, this choice of RPs introduces projection errors.
    The more general the hull type is, the fewer RPs cover the dataset.
  }
  \label{fig:hull-types}
\end{figure*}

We illustrate the differences between these hulls in Figure~\ref{fig:hull-types}, where each successive hull type expands the space in which the original data points can be expressed using the RPs without introducing a projection error.

\subsubsection{Greedy Hull Clustering}

While the concept of using a convex or conic hull to minimize the projection error is natural, identifying such a hull in practice is nontrivial.
Specifically, selecting a small subset of representative periods whose hull contains all the original periods is a computationally hard problem.
This is related to well-known non-negative matrix factorization problem, which is NP-hard \cite{vavasis2009}.

To address this challenge efficiently, we adapt the greedy conic hull algorithm \cite[Algorithm~1]{Buskirk} for all three hull types.
The algorithm identifies points in the conic hull via a reduction to the convex hull with the gnomonic projection (i.e., a projection of rays onto points in a plane);
omitting the projection yields a direct convex-hull method.

We additionally improve the algorithm by employing caching to avoid redundant distance computations during the greedy selection of the hull points.
Specifically, once a point’s distance to the current convex hull has been computed, that value is stored in a cache.
Since the convex hull can only grow as new points are added, the distance from any given point to the hull cannot increase over time.
This property allows the algorithm to safely reuse previously computed distances unless there is strong evidence that a point’s position relative to the hull has changed significantly---namely, if its distance to the newly added point is smaller than its cached hull distance.
By using this check, the algorithm reduces the number of expensive projection operations, thereby improving overall efficiency.
Although this criterion is a heuristic, we nevertheless observed that in practice it offers a substantial speedup.
The resulting method is presented in Algorithm~\ref{alg:GreedyHull}.

\begin{algorithm}[htb]
\caption{%
Greedy Convex Hull Clustering\label{alg:GreedyHull}
}
\KwData{%
  clustering matrix \(\mathbf{C}\),
  distance function \(\dist\),\\
  hull projection operator \(\operatorname{Hull}\),
  initial RP set \(\R\),\\
  number of representative periods \(\numRPs\).
}
\KwResult{%
  a set of RPs \(\R \subseteq \D\).
}
\If{\(\R\) is empty}{
compute the mean \(\bar{\mathbf{c}} \gets \frac{1}{D} \sum_{d \in \D} \mathbf{c}_{d}\)\;
select the point furthest from the mean, \(r \gets \argmax_{d \in \D} \dist(\mathbf{c}_{d}, \bar{\mathbf{c}})\)\;
initialize the RP index set, \(\R \gets \{r\}\)\;
}
initialize an empty distance cache \(\mathfrak{C} \colon \D \to \Reals_+\)\;
\While{\(\nR < \numRPs\)}{
  let \(\newr\) be the most recently added index in \(\R\)\;
  \(\maxdist \gets -\infty\)\;

  \ForEach{\(d \in \D \setminus \R\)}{
    \eIf{\(d \in \mathfrak{C}\) \textbf{and}
          \(\dist\!\bigl(\mathbf{c}_d,\mathbf{c}_{\newr}\bigr) \geq \mathfrak{C}(d)\)}{
      %--- cached value is still a valid upper bound
      \(\curdist \gets \mathfrak{C}(d)\)\;
    }{
      %--- recompute expensive projection
      \(\curdist \gets 
        \dist \bigl(\mathbf{c}_d,
        \operatorname{Hull}\{\mathbf{c}_r \mid r\in\R\}\bigr)\)\label{line:prject-onto-hull}\;
        \(\mathfrak{C}(d) \gets \curdist\)\;
    }
    \If{\(\curdist > \maxdist\)}{
      \(\maxdist \gets \curdist,\)\quad\(\maxr \gets d\)\;
    }
  }

  \(\R \gets \R \cup \{\maxr\}\)\;
}

\Return{\(\R\).}
\end{algorithm}

Depending on the hull type, the full clustering procedure is:
\begin{itemize}
  \item \emph{Convex hull}: Algorithm~\ref{alg:GreedyHull} is applied
        directly to \(\mathbf{C}\).
  \item \emph{Convex hull with null}:
    \begin{enumerate}
      \item add a null vector \(\mathbf{0}\) as a column to \(\mathbf{C}\),
      \item run Algorithm~\ref{alg:GreedyHull} with \(\nR + 1\) representatives and initial set of RPs \(\R = \{\mathbf{0}\}\), then
      \item drop the null vector \(\mathbf{0}\) from both \(\mathbf{C}\) and \(\R\).
    \end{enumerate}
  \item \emph{Conic hull}: use the gnomonic projection \cite{Buskirk}:
    \begin{enumerate}
      \item compute the column mean \(\bar{\mathbf{c}}\) of the data and normalize it to obtain a reference direction \(\mathbf{q}\),
      \item for each column \(\mathbf{c}_{d}\) of the clustering matrix \(\mathbf{C}\), scale it by \(1 / \langle \mathbf{c}_{d}, \mathbf{q} \rangle\) so that all points lie on the hyperplane \(\mathcal{H} = \bigl\{\mathbf{x} \bigm| \langle \mathbf{x}, \mathbf{q}\rangle = 1\bigr\}\),
      \item apply Algorithm~\ref{alg:GreedyHull} to the scaled clustering matrix,
      \item return the corresponding unscaled RPs.
    \end{enumerate}
\end{itemize}

The greedy procedure ensures that each added representative improves coverage of the dataset, and it can be terminated once a fixed number of representatives is reached.

At each iteration, the method aims to select the data point that most reduces the projection error of the remaining points onto the current hull.
This approach allows us to build up the RP set incrementally while maintaining a good approximation of the original data within the chosen hull type.

Computing the distance from a data point \(\mathbf{c}_d\) to the current convex hull \(\operatorname{Hull}\{\mathbf{c}_r \mid r \in \R\}\), as required in line~\ref{line:prject-onto-hull} of Algorithm~\ref{alg:GreedyHull}, is a nontrivial task.
We formulate this as a projection problem: given a point \(\mathbf{c}_d\), we seek the closest point in the convex hull in terms of Euclidean distance \(\dist(\mathbf{R}\mathbf{w}_{\!d}, \mathbf{c}_{d}) = \lVert\mathbf{R}\mathbf{w}_{\!d} - \mathbf{c}_{d}\rVert_2\).
Since any point in the hull of the columns of \(\mathbf{R}\) can be expressed as a convex combination of those columns using weights \(\mathbf{w}_{\!d}\), that is, \(\mathbf{R}\mathbf{w}_{\!d}\), this projection reduces to the following constrained optimization problem:
\begin{mini!}
{}
{\tfrac{1}{2} \lVert\mathbf{R}\mathbf{w}_{\!d} - \mathbf{c}_{d}\rVert_2^2\label{obj:pgd-clustering}}{\label{eq:weight-fitting-problem}}{}
\addConstraint{\mathbf{w}_{\!d}}{\in \mathbb{W}^{\text{conv}}.\label{eq:weight-fitting-problem-const}}
\end{mini!}
Without the constraint \eqref{eq:weight-fitting-problem-const}, the problem has a well-known solution \(\mathbf{R}^\dagger \mathbf{c}_d\) where \(\mathbf{R}^\dagger\) is the Moore--Penrose pseudoinverse of \(\mathbf{R}\).
This might result in negative weights, and even projecting back onto the subspace of acceptable weights \(\proj_{\mathbb{W}^{\text{conv}}} (\mathbf{R}^\dagger \mathbf{c}_d)\) yields a good feasible but not necessarily optimal solution.

To address this issue, we use projected gradient descent (PGD, described in Algorithm~\ref{alg:PGD}) to further fine-tune the initial weight guess.
Note that replacing \(\lVert\placeholder\rVert_2\) with \(\tfrac{1}{2}\|\cdot\|^2_2\) does not alter the location of the minimizer, since the function \(f(x) = \tfrac{1}{2}x^2\) is strictly increasing on \(\Reals_+\).
Thus, minimizing the squared distance is equivalent to minimizing the distance itself, while offering a simpler  subgradient \(g\) of the objective:
\begin{equation}
  g(\mathbf{w}_{\!d}) = \nabla_{\!\mathbf{w}_{\!d}} \tfrac{1}{2}\lVert\mathbf{R}\mathbf{w}_{\!d} - \mathbf{c}_{d}\rVert_2^2 = \mathbf{R}^\top(\mathbf{R}\mathbf{w}_{\!d} - \mathbf{c}_{d}).
  \label{eq:subgrad}
\end{equation}
This subgradient is smooth, making the projection step more amenable to efficient optimization using PGD, which we employ to compute the distance to the hull in line~\ref{line:prject-onto-hull} of Algorithm~\ref{alg:GreedyHull} as follows:
\begin{align*}
  \curdist &\gets \dist\bigl(\mathbf{c}_d, \operatorname{Hull}\{\mathbf{c}_r \mid r\in\R\}\bigr)
  \approx
  \dist (\mathbf{c}_d, \mathbf{R}\mathbf{w}_{\!d}),\;\text{where}\\
  \mathbf{w}_{\!d}
  &\gets \pgd\bigl(\proj_{\mathbb{W}^{\text{conv}}} (\mathbf{R}^\dagger \mathbf{c}_d) \bigm| g, \proj_{\mathbb{W}^{\text{conv}}}, N, \varepsilon, \alpha\bigr).
\end{align*}
We use \(\proj_{\mathbb{W}^{\text{conv}}} (\mathbf{R}^\dagger \mathbf{c}_d)\) as the initial weight guess and \eqref{eq:subgrad} as the subgradient function \(g\).
For the projection operator \(\proj_{\mathbb{W}^{\text{conv}}}\),
we utilize the fast algorithm \cite[Figure 2]{Condat2016} for projection onto the simplex \(\bigtriangleup_n = \mathbb{W}^{\text{conv}}\).
For the implementation details, the reader is referred to the original paper.

\begin{algorithm}[htb]
\caption{%
Projected Gradient Descent\\%
\(\pgd(\mathbf{x} \mid g, \proj_{\mathbb{X}}, N, \varepsilon, \alpha)\)\label{alg:PGD}
}
\KwData{%
  vector \(\mathbf{x} \in \Reals^n\),
  subgradient \(g: \Reals^n \to \Reals^n\),\\
  projection operator \(\proj_{\mathbb{X}}\),
  maximum number of iterations \(N\),
  tolerance \(\varepsilon\),
  learning rate \(\alpha\).
}
\KwResult{%
  an approximate projection of \(\mathbf{x}\) onto \(\mathbb{X}\).
}
project the initial point \(\mathbf{x}\) onto space \(\mathbb{X}\),  \(\mathbf{x} \gets \proj_{\mathbb{X}}\mathbf{x}\)\;
\For{\(N\) iterations}{
  compute the subgradient \(\mathbf{g}\) at point \(\mathbf{x}\), \(\mathbf{g} \gets g(\mathbf{x})\)\;
  remember the current point, \(\mathbf{x}_{\text{prev}} \gets \mathbf{x}\)\;
  descent against the subgradient, \(\mathbf{x} \gets \mathbf{x} - \alpha \mathbf{g}\)\;
  project point \(\mathbf{x}\) onto space \(\mathbb{X}\), \(\mathbf{x} \gets \proj_{\mathbb{X}}\mathbf{x}\)\;
  \If{the point didn't move, \(\lVert\mathbf{x}_{\text{prev}} - \mathbf{x}\rVert_\infty \leq \varepsilon / N\)}{
    \Break\;
  }
}
\Return{the resulting point \(\mathbf{x}\).}
\end{algorithm}

\subsection{Finding the Blended Weights}

After hull clustering is complete, we obtain a matrix \(\mathbf{R}\) of the RP data.
The next step is to find the weights \(\pWeightRepr{d,r}\) that best approximate the original data for each period \(d\).

This is the weight-fitting problem \eqref{eq:weight-fitting-problem} but with the weights potentially restricted to a different domain \(\mathbb{W}^{\placeholder}\). It can be solved with Algorithm~\ref{alg:PGD} using different projection operators \(\proj_{\mathbb{W}^{\placeholder}}\) depending on the weight type.
\begin{itemize}
  \item \emph{Dirac weights} \(\mathbb{W}^{\text{Dirac}}\): project by setting the largest coordinate of \(\mathbf{x}\) to 1 and all others to 0. In practice, cluster assignments directly provide this projection, so explicit projection is often unnecessary.

\item \emph{Convex weights} \(\mathbb{W}^{\text{conv}}\): apply the aforementioned fast projection algorithm \cite[Figure 2]{Condat2016}.

\item \emph{Conic sub-unit weights} $\mathbb{W}^{\leq 1}$: use the convex hull method with the null point trick as described previously.

\item \emph{Conic weights} \(\mathbb{W}\): project by thresholding negative values to zero,
\(
\proj_{\mathbb{W}}(\mathbf{x}) = \mathbf{x}^+ = \max(\mathbf{x}, \mathbf{0}).
\)
\end{itemize}

For the initial weight guess, we use either \(\proj_{\mathbb{W}^{\text{conv}}} (\mathbf{R}^\dagger \mathbf{c}_d)\) or the Dirac weights returned by the clustering algorithms, whichever results in a smaller projection error.
\section{Case Studies}

We test the proposed methodology in two case studies.\footnote{The input data and implementation for both case studies are available online at \url{https://github.com/greg-neustroev/blended-rep-periods}.}
The first case study is a GEP problem based on information from the TYNDP 2022 \cite{tyndp}.
This stylized case study allows for investments in renewable energy from greenfield setting.
For the sake of simplicity, this case study includes ramping constraints but no storage and conversion technologies.

The second case study is also based on information from the TYNDP, but from the 2024 study;
in this case, we include hydrogen demand for each node, which can be met through electrolyzers that convert electricity into hydrogen or through steam methane reforming.
Unlike the GEP case study, this analysis excludes ramping constraints and investment decisions, with all capacities predetermined. It includes both short-term and seasonal storage options, such as batteries and pumped hydro storage, and models hydropower plants, reservoirs, and inflows in each country. This approach highlights the importance of seasonal behavior, particularly in the Nordic countries.

\subsection{Experiments Setup}

For each case study, we design experiments to evaluate the performance of different clustering methods and weight types.
\begin{itemize}
\item \emph{Weight types}: we compare standard \textit{Dirac weights} with the proposed \textit{conic}, \textit{conic sub-unit}, and \textit{convex} weights.
\item \emph{Clustering methods}: in addition to the common approaches (\textit{\(k\)-means}, \textit{\(k\)-medoids}), we consider \textit{hull} methods specific to each weight type (convex for Dirac).
\item \emph{Number of representatives}: we vary the number of RPs to analyze its impact on each method’s performance.
\end{itemize}

Each experiment is run with five different random seeds to account for various sources of randomness.

\subsection{Regret as a Measure of Solution Quality}

RP methods approximate full-resolution ESOMs---e.g., modeling a full year hourly---but it can be hard to judge how good these approximations are.
Rather than just comparing objective values, we use \emph{relative regret}, which measures the extra cost (in percent) incurred when decisions are made using the approximate model instead of the full one.
In other words, regret shows how much worse the solution is because of the simplification, with lower regret values being better.

In the GEP case study, we compute the regret of the investment decisions, that is, the investment variables from the approximate model are fixed, and the full ESOM is rerun to calculate the actual total cost (i.e., investment and operation).
These costs are compared against the costs of the full original ESOM.

In the P2X case, seasonal storage levels (states of charge at the end of every day) are fixed, and rerunning the full problem gives the operational cost; regret is again the difference from the benchmark.

\subsection{Results for Generation Expansion Planning}

We run the GEP case studies for 5, 10, 20, 40, and 80 RPs. 

\begin{figure*}[tbh]
\centering
\includegraphics[width=\textwidth]{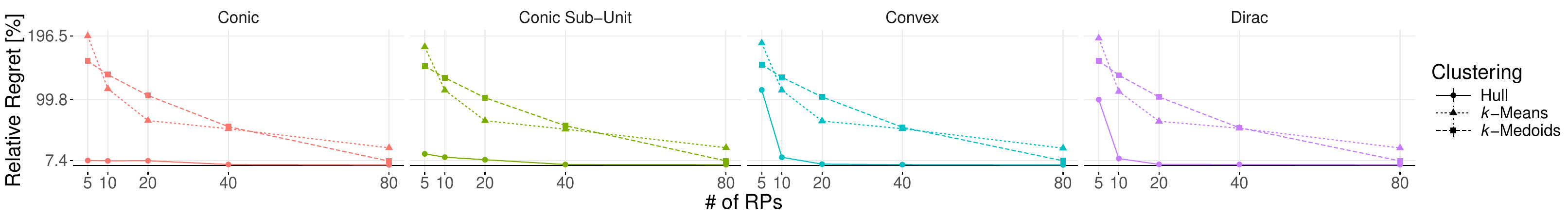}
\caption{
    Relative regret per clustering method in the GEP case study for four types of weights. Hull clustering outperforms \(k\)-means and \(k\)-medoids.
  }
  \label{fig:gep-regret}
\end{figure*}

\begin{figure}[tbh]
\centering
\includegraphics[width=0.488\textwidth]{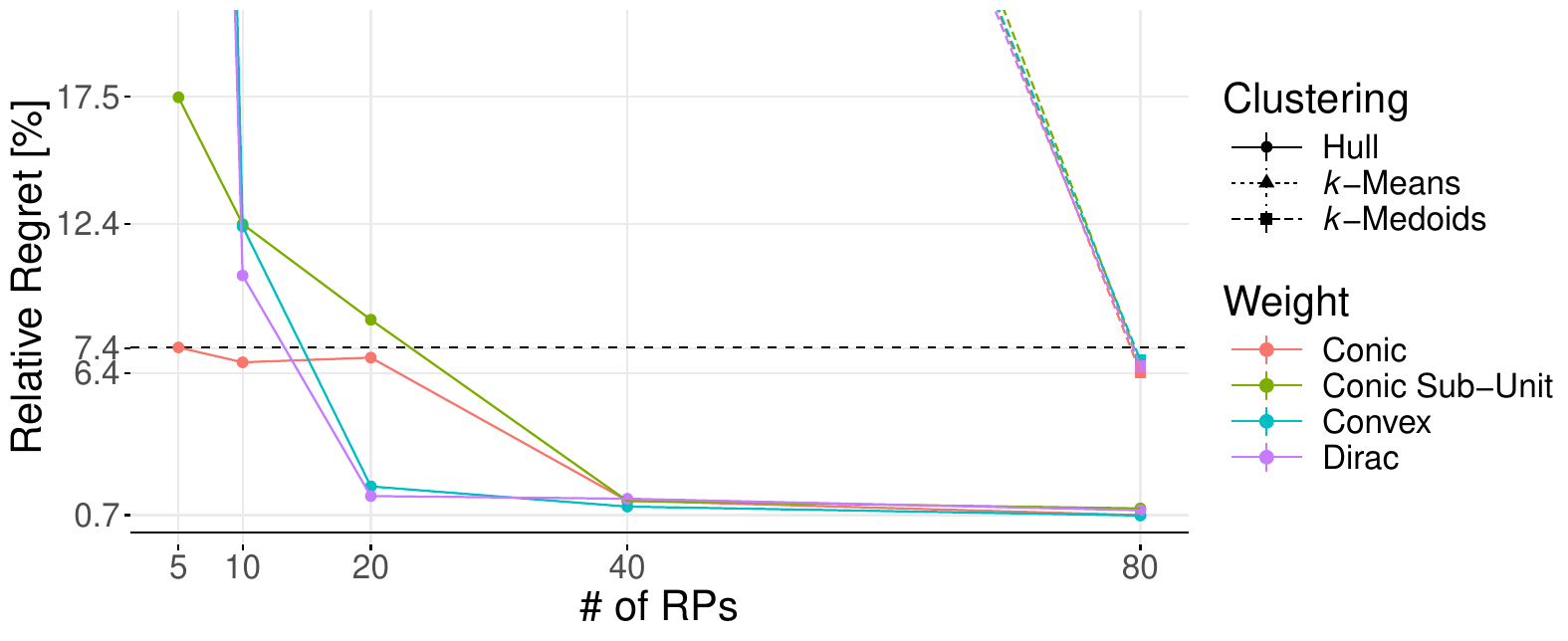}
\caption{
    Relative regret in the GEP case study for best-performing methods.
    The dashed line shows the relative regret of 7.4\% with just 5 RPs, which is only achieved for non-hull methods for 80+ RPs.
  }
  \label{fig:gep-regret-hull}
\end{figure}

Figure~\ref{fig:gep-regret} shows the regret of all of the experiments grouped by the weight type.
In case of both \(k\)-means and \(k\)-medoids, the weight type has almost no effect on the performance.
In each case, the hull method outperforms both conventional clustering methods.
Moreover, with just 5 RPs, conic clustering achieved relative regret of 7.4\%, a level not achieved by any of the \(k\)-means results (lowest relative regret 26.1\%) and matched only by \(k\)-medoids with 80 RPs (6.4\%), as seen in Figure~\ref{fig:gep-regret-hull}.

\begin{figure}[tbh]
\centering
\includegraphics[width=0.488\textwidth]{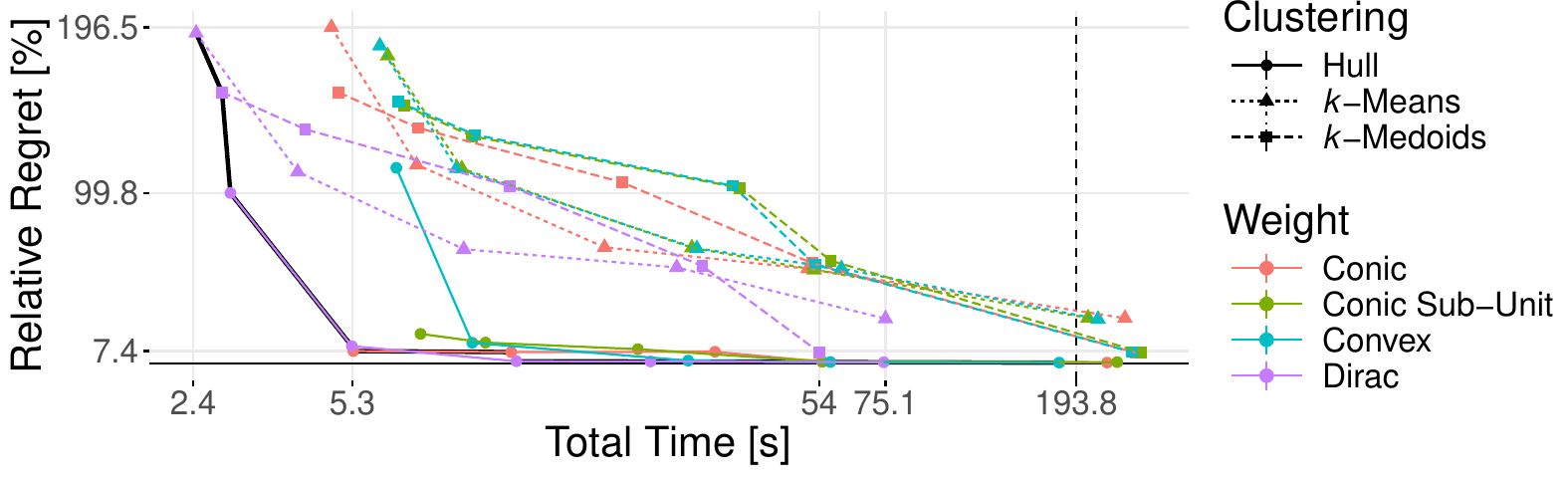}
\caption{
    Relative regret versus (logarithmic) total time for GEP. When relative regret is under 100\%, the Pareto front (thick black line) consists exclusively of hull methods.
    The coordinates of break points in the Pareto front are used as labeled tick marks on the axes.
    The dashed line shows mean time to solve the full problem.
  }
  \label{fig:gep-time-regret}
\end{figure}

Figure~\ref{fig:gep-time-regret} shows the total time, that is, the time required to read, preprocess, and cluster the data, fit the weights, formulate the model and solve it.
Conic hull with 5 RPs took \(5.3 \pm 0.4\) s, compared to \(54.0 \pm 0.6\) s for 80-RP \(k\)-means;
it was, therefore, ten times faster and achieved 3.5 times lower regret.
\(k\)-medoids with 80 RPs took \(75.1 \pm 1.9\) s, being 14 times slower than conic hull method to achieve similar regret levels.

Blended weight methods, however, scale worse than Dirac methods: at 80 RPs they were slower than solving the full problem directly, as weight fitting dominates runtime.
In contrast, Dirac weights with twice as many RPs often matched blended performance, highlighting the trade-off between weight sophistication and the number of RPs.

In summary, in GEP hull clustering---especially with conic weights---provided superior accuracy and efficiency, dominating both \(k\)-means and \(k\)-medoids.

\subsection{Results for Power-to-X Dispatch with Seasonal Storage}

\begin{figure*}[tbh]
\centering
\includegraphics[width=\textwidth]{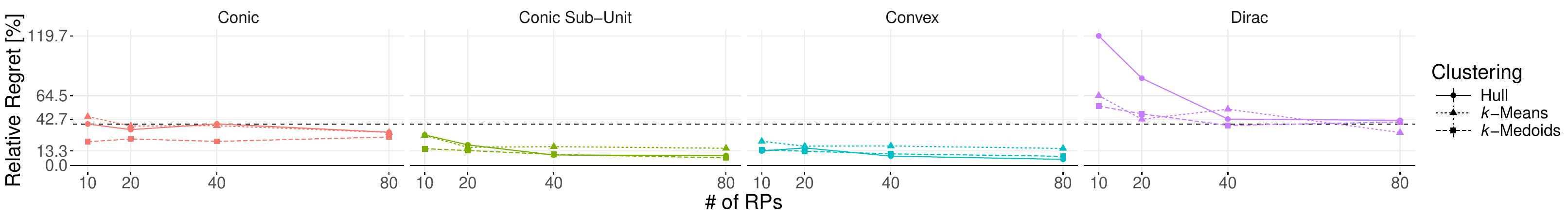}
\caption{
    Relative regret per clustering method in the P2X case study for four types of weights. Blended weights outperform Dirac weights. The worst performance of hull clustering with blended RPs is shown in dashed black line.
  }
  \label{fig:p2x-regret}
\end{figure*}

\begin{figure*}[bth]
\centering
\includegraphics[width=\textwidth]{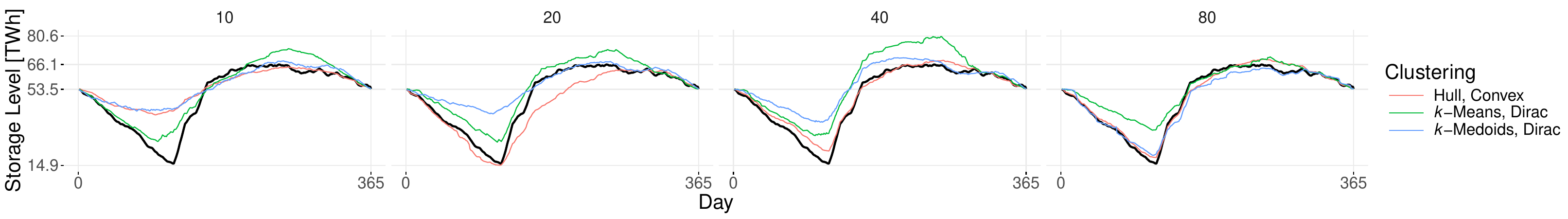}
\caption{
    Inter-period hydro reservoir storage levels for the largest reservoir (Norway) for different numbers of RPs. Thick black line shows the full solution.
  }
  \label{fig:p2x-reservoir}
\end{figure*}

As the P2X case study is more complex, 5 RPs were insufficient; therefore, experiments began with 10 RPs.

\begin{figure}[tbh]
\centering
\includegraphics[width=0.488\textwidth]{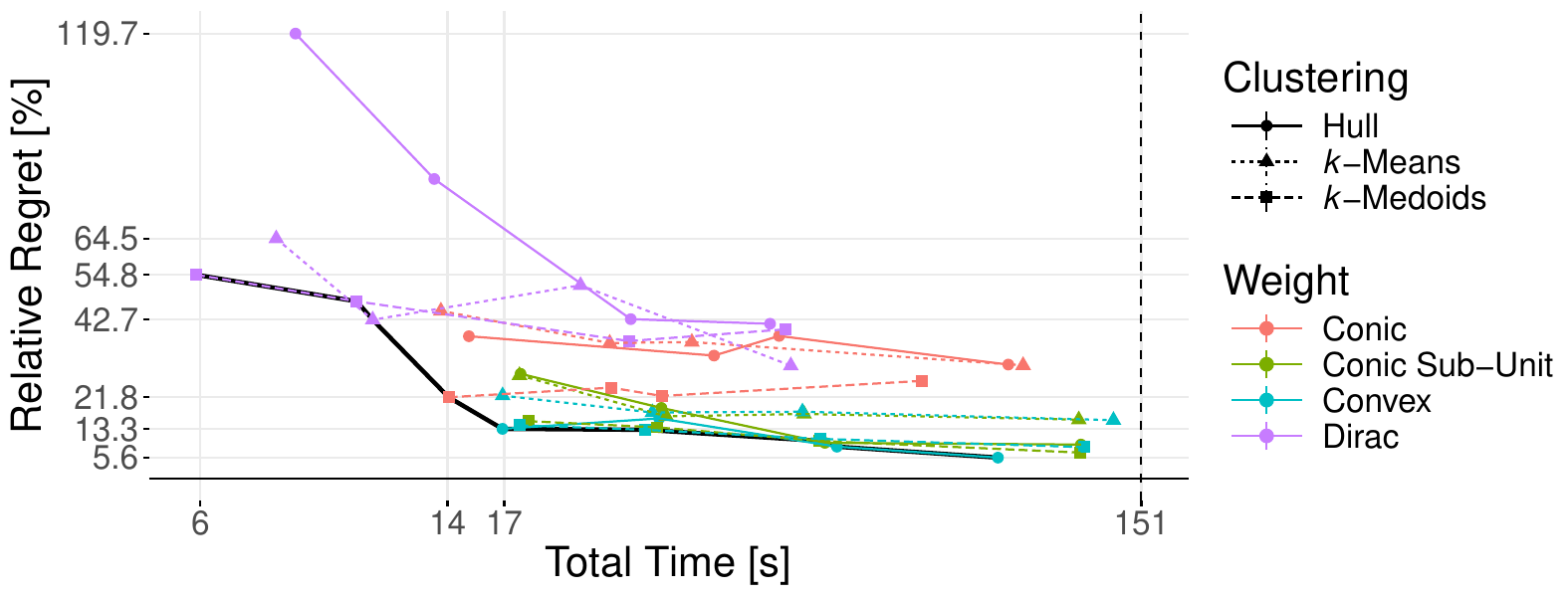}
\caption{
    Relative regret versus (logarithmic-scale) total time in the P2X case study. The Pareto front is shown in thick black line.
    The coordinates of break points in the Pareto front are used as labeled tick marks on the axes.
    The dashed line shows mean time to solve the full problem.
  }
  \label{fig:p2x-time-regret}
\end{figure}

Figure~\ref{fig:p2x-regret} reports relative regret for all methods.
Unlike in the GEP case, weight type plays a decisive role: blended weights outperform Dirac weights, regardless of clustering method.
The runtime–regret trade-off (Figure~\ref{fig:p2x-time-regret}) highlights two parts of the Pareto front: convex-weight hull clustering, which achieves regret below 15\% with a runtime of 17 s, and Dirac \(k\)-medoids, which runs faster but with more than 50\% regret.

Blended weights also capture reservoir dynamics better. Figure~\ref{fig:p2x-reservoir} shows seasonal storage levels for Norway, the largest reservoir in the system.
From 40 RPs onward, hull clustering closely tracks the full model, whereas \(k\)-means diverges.
With just 10 or 20 RPs, convex hull clustering yields near-optimal solutions (13.3\% and 16.1\% regret) compared to much higher regret levels for \(k\)-means (64.5\% and 42.7\%) and \(k\)-medoids (54.8\% and 47.6\%) with the same numbers of RPs.

Overall, in P2X the choice of weight type is more important, with convex hull methods striking the best balance between regret, runtime, and reservoir fidelity.
\section{Discussion}

The results highlight several benefits of using RP methods with advanced weighting and clustering strategies.
Hull clustering with blended weights consistently delivers lower regret, faster runtimes, and more faithful approximations.
These improvements are especially valuable when solving large-scale, high-resolution ESOMs, where runtime and memory requirements would otherwise be prohibitive.

At the same time, the case studies reveal that performance is strongly problem-dependent.
In GEP with ramping, the main challenge is capturing chronological variability;
what matters most is selecting RPs that preserve extreme ramps.
Clustering type therefore dominates performance: hull methods excel because they select diverse boundary points that best capture these extremes.
By contrast, in P2X, the binding constraints are less about ramps and more about seasonal storage balances.
Once RPs are chosen, weights determine how they are ``stitched back together'' into a synthetic year.
Here, the choice of weights is critical.
Even though clustering type is more important for GEP and weight type for P2X, the combination of hull clustering and blended weights performs well in both cases;
while not necessarily optimal, using them together appears to be important, as they exploit different problem structures.

The experiments highlight potential pitfalls.
The problem-dependence means that applying a method blindly across models may lead to poor results.
Moreover, blended weights, while powerful, introduce additional computational overhead; at the same time, we do not expect this to be an issue for large-scale problems, since the time to solve the optimization problem will very likely dominate the clustering time.
For very large numbers of RPs, weight fitting becomes the bottleneck, sometimes exceeding the runtime of the full problem.
This makes parameter tuning, such as the learning rate and step length in the gradient descent algorithm used for weight fitting, critical to achieving convergence without excessive iterations, calling for further research.
\section{Conclusion}
Existing clustering-based approaches using representative periods (RPs) to solve large ESOMs can be captured in a general framework.
Within this framework, we introduced hull clustering with blended RPs, a novel approach to selecting and combining extreme and typical time series to efficiently reduce temporal data in ESOMs.
By focusing on extreme points and blending them, the method captures critical system conditions more accurately than conventional clustering techniques, while requiring fewer RPs.
Through two case studies---a stylized GEP problem and a more comprehensive P2X dispatch problem---we demonstrated that the approach achieves high solution quality with substantially reduced computational burden.
The results highlight its ability to preserve important behavior patterns, improving fidelity, and achieving lower regret.
While implementing blended RPs introduces new challenges in the clustering methodology, the approach offers a promising avenue for scaling ESOMs to longer horizons and more detailed models.
Future work could extend the method to a more general formulation, test it on more diverse case studies, and integrate it directly into existing ESOM tools.

\section*{Acknowledgments}
This publication is part of the project NextGenOpt with project number ESI.2019.008, which is financed by the Dutch Research Council (NWO) and supported by eScienceCenter under project number NLeSC C 21.0226. In addition, this research received partial funding from the European Climate, Infrastructure and Environment Executive Agency under the European Union’s HORIZON Research and Innovation Actions under grant agreement \textnumero\ 101095998.

\bibliographystyle{IEEEtran}
\bibliography{literature.bib}

\end{document}